\documentclass[12pt,a4paper]{article}
\usepackage{amssymb,amsmath}
\usepackage{amsfonts}
\usepackage{amsthm}
\usepackage{a4}
\usepackage{tikz}
\usepackage{verbatim}
\usepackage{caption}
\usepackage{float}
\newcommand{\para}{\par\vspace{.25cm}}
\newtheorem{defi}{Definition}
\newtheorem{prop}{Proposition}
\newtheorem*{theorem*}{Theorem}
\newtheorem{theorem}{Theorem}
\newtheorem{lemma}{Lemma}
\newtheorem{cor}{Corollary}

\newcommand{\Q}{\mathbb{Q}}

\newcommand{\Z}{\mathbb{Z}}

\newcommand{\C}{\operatorname{Cen}}

\begin{document}
\baselineskip 18pt \title{{\bf \ Central units of integral group rings of monomial groups} {} { }}
\author{ Gurmeet K. Bakshi   {\footnote {Research supported by Science and Engineering Research Board (SERB), DST, Govt. of India under the scheme Mathematical Research Impact Centric Support (sanction order no MTR/2019/001342) is gratefully acknowledged.}}\\ {\em \small Centre for Advanced Study in
Mathematics,}\\
{\em \small Panjab University, Chandigarh 160014, India}\\{\em
\small email: gkbakshi@pu.ac.in}\vspace{.2cm}\\ and \vspace{.2cm}\\ Gurleen Kaur {\footnote {Research supported by Science and Engineering Research Board (SERB), DST, Govt. of India (National Post-Doctoral Fellowship Reference No. PDF/2020/000343) is gratefully acknowledged.} \footnote{Corresponding author}}\\{\em \small Department of Mathematical Sciences,}\\ {\em \small Indian Institute of Science Education and Research Mohali,}\\{\em \small Sector 81, Mohali (Punjab) 140306, India}\\{\em
\small email: gurleen@iisermohali.ac.in}}
\date{}
{\maketitle}
\begin{abstract}\noindent In this paper, it is proved that the group generated by Bass units contains a subgroup of finite index in the group of central units $\mathcal{Z}(\mathcal{U}(\mathbb{Z}G))$ of the integral group ring $\mathbb{Z}G$ for a subgroup closed monomial group $G$ with the property that every cyclic subgroup of order not a divisor of $4$ or $6$ is subnormal in  $G$. If $G$ is a generalized strongly monomial group, then it is shown that the group generated by generalized Bass units contains a subgroup of finite index in $\mathcal{Z}(\mathcal{U}(\mathbb{Z}G))$. Furthermore, for a generalized strongly monomial group $G$, the rank of  $\mathcal{Z}(\mathcal{U}(\mathbb{Z}G))$ is determined. The formula so obtained is in terms of generalized strong Shoda pairs of $G$.\end{abstract}\section{Introduction}Throughout this paper, $G$ denotes a finite group. Let  $\mathcal{U}(\mathbb{Z}G)$ be the unit group of the integral group ring $\Z G$ and let  $\mathcal{Z}(\mathcal{U}(\mathbb{Z}G))$  be its center. It is well known that $\mathcal{Z}(\mathcal{U}(\mathbb{Z}G))$ is a finitely generated abelian group.
In this paper, we are concerned with the rank of  $\mathcal{Z}(\mathcal{U}(\mathbb{Z}G))$ and large subgroups of $\mathcal{Z}(\mathcal{U}(\mathbb{Z}G))$, i.e., subgroups of finite index in $\mathcal{Z}(\mathcal{U}(\mathbb{Z}G))$.  Both of these problems have been the center of attraction for several decades starting with the work of Higman \cite{Higman} and require a deep understanding of the structure of $G$ and that of the rational group algebra $\Q G$. Although a lot of information is known, yet it is far from being completely understood. We refer to the various surveys including dedicated books on the topic \cite{Jespers1,JdR,Milies1,MS,Sehgal, Sehgal1} for its history. Bass units (introduced by Bass \cite{Bass}) are generic  construction of units in $\Z G$ which are analogous  to cyclotomic units.  Bass and Milnor \cite{Bass} proved that if $G$ is abelian, then the group generated by Bass units contains a large subgroup of  $\mathcal{Z}(\mathcal{U}(\mathbb{Z}G))$.
Except some special cases, Bass units are of infinite order and are not generally central units. However, one can construct central units of $\Z G$ by taking product of conjugates of Bass units with the help of suitable subnormal series in $G$, if available. This type of construction is seen in the work of Jespers, Parmenter and Sehgal in \cite{JPS} where it is proved that if $G$ is a nilpotent group, then the group generated by Bass units contains a large subgroup of $\mathcal{Z}(\mathcal{U}(\mathbb{Z}G)).$ This result was extended by Jespers, Olteanu, del R{\'{\i}}o and Van Gelder in \cite{JOdRV} to abelian-by-supersolvable groups $G$ with the property that every cyclic subgroup of order not a divisor of 4 or 6 is subnormal in $G$.  In a recent work \cite{BK2}, we gave a further extension of this result to the class $\mathcal{C}$ which consists of all finite groups whose each subquotient is either abelian or contains a non central abelian normal subgroup, but with a constraint on a complete and irredundant set of Shoda pairs. In this paper, we will show that the imposed constraint  is no longer needed and in fact something  better holds. In the arguments given in all the papers \cite{BK2,JPS,JOdRV} one strongly require the underlying class of groups to be closed under subgroups because the idea is to proceed by induction. Thus, if we confine ourselves to monomial groups where there is a strong interplay between the subgroup structure of $G$  and the structure of $\Q G$,  the maximum  potential where the ideas in \cite{BK2,JPS,JOdRV} can be generalized is the class of subgroup closed monomial groups and in this paper, we prove the following:\begin{quote} If $G$ is a subgroup closed monomial group such that every cyclic subgroup of order not a divisor of $4$ or $6$ is subnormal in $G$, then  the group generated by Bass units  of $\mathbb{Z}G$ contains a large subgroup of $\mathcal{Z}(\mathcal{U}(\mathbb{Z}G))$. \end{quote}
It  may be pointed out that  nilpotent groups  $\subsetneq$ abelian-by-supersolvable groups \linebreak $\subsetneq$   $\mathcal{C}$  $\subsetneq$   subgroup closed monomial groups   (see Proposition 1 of \cite{BK1}).  \par It is quite natural to ask if one can  construct units in $\Z G$ from those of $\Z G/N$, i.e., the integral group ring of a quotient group of $G$.  In this connection, one can show that  if $b$ is a unit in  $\Z G/N$, then some power of $b$ has a pre-image in $\Z G$ which is a unit. Utilizing  this idea,  in \cite{JOdRV}, generalized Bass units are defined by taking $b$ to be Bass units in $ \Z G/N$. There it is proved that the group generated by generalized Bass units of $\Z G$ contains a large subgroup of  $\mathcal{Z}(\mathcal{U}(\mathbb{Z}G))$, if $G$ is a strongly monomial group. This result was extended in \cite{BK2} to generalized strongly monomial groups with a constraint on a complete and irredundant set of Shoda pairs.  Here, we will evade that constraint and prove the following: \begin{quote} If $G$ is a generalized strongly monomial group, then the group generated by generalized Bass units of $\mathbb{Z}G$ contains a large subgroup of $\mathcal{Z}(\mathcal{U}(\mathbb{Z}G))$. \end{quote} One of the important steps to compute large subgroups of $\mathcal{Z}(\mathcal{U}(\mathbb{Z}G))$  is to know its rank.  In this connection, it is known from the works of  \cite{Ferraz} and  \cite{RS} that  the rank of $\mathcal{Z}(\mathcal{U}(\mathbb{Z}G))$  for any finite group $G$  is the difference between the number of  simple components of the real group algebra $\mathbb{R} G$ and that of $\Q G$. Furthermore, the number of simple components of  $\mathbb{R} G$ (resp. $\Q G$) coincides with the number of real (resp. rational) conjugacy classes of $G $ (see Theorem 42.8 of \cite{Curtis}). For metacyclic groups, Ferraz and Simon in \cite{FS} gave a precise formula  for the rank in terms of the order of $G$. For strongly monomial groups, a description of the rank in terms of strong Shoda pairs of $G$ was given by  Jespers, Olteanu, del R{\'{\i}}o and Van Gelder in \cite{JOdRVG1}. In this paper, we will  extend this work  and  provide a  precise count of the rank of  $\mathcal{Z}(\mathcal{U}(\mathbb{Z}G))$ when $G$ is generalized strongly monomial. It is pertinent to mention that the class of generalized strongly monomial groups  include all well known classes of monomial groups like abelian-by-supersolvable, strongly monomial, $\mathcal{C}$, subnormally monomial, supermonomial, monomial groups with Sylow towers and many more. For details on the vastness of the class of generalized strongly monomial groups, we refer to \cite{BK4}. \section{ Background on Shoda pair theory} \noindent In this section, we quickly recall the basic definitions of Shoda pair theory and some important results needed for our purpose.  Shoda pairs and strong Shoda pairs were introduced by Olivieri, del R{\'{\i}}o and Sim{\'o}n in  \cite{OdRS04}.  If $\chi$ is a monomial character of $G$, one would like to know a subgroup $H$ of $G$ and $K \unlhd H$ such that $\chi$ can be induced from a linear character of $H$ with kernel $K$. One is interested in determining pairs of subgroups $(H,K)$ of $G$ such that a linear character of $H$ with kernel $K$ when induced to $G$ is irreducible. This is precisely the work of Shoda (see  \cite{JdR}, Corollary 3.2.3) and in his honor Olivieri, del R{\'{\i}}o and Sim{\'o}n \cite{OdRS04} termed such pairs as Shoda pairs of $G$. More precisely, a \textit{Shoda pair} (\cite{OdRS04}, Definition 1.4) of $G$ is a pair $(H, K)$  of subgroups of $G$ satisfying the following:\begin{itemize}
	\item [(i)] $K \unlhd H$,  $H/K$ is cyclic;
	\item [(ii)] if $g \in G$ and $[H, g] \cap H \subseteq K,$ then $ g \in H.$\end{itemize}
 For $K\unlhd H\leq G$, define:$$\widehat{H}:=\frac{1}{|H|}\displaystyle\sum_{h \in H}h,$$ $$\varepsilon(H,K):=\left\{\begin{array}{ll}\widehat{K}, & \hbox{$H=K$;} \\\prod(\widehat{K}-\widehat{L}), & \hbox{otherwise,}\end{array}\right.$$ where $L$ runs over all the minimal normal subgroups of $H$ containing $K$ properly, and $$e(G,H,K):= {\rm~the~sum~of~all~the~distinct~}G{\rm {\tiny{\operatorname{-}}} conjugates~of~}\varepsilon(H,K).$$ \para \noindent  If $(H, K)$ is a Shoda pair of $G$ and $\lambda$ is a linear character of $H$ with kernel $K$,  denote by  $e_{\mathbb{Q}}(\lambda^G)$, the primitive central idempotent of $\Q G$ corresponding to the irreducible character $\lambda^G$.  It is proved in  Theorem 2.1 of \cite{OdRS04} that $e_{\mathbb{Q}}(\lambda^G)$, which is equal to $\frac{1}{|H|}\sum_{\sigma \in \operatorname{Gal}(\mathbb{Q}(\lambda^{G})/\mathbb{Q})}\sum_{g \in G}\sigma
(\lambda^{G}(g))g^{-1}$, is a rational multiple (unique) of $e(G, H, K)$.  So the simple component $\Q G e_{\mathbb{Q}}(\lambda^G)$  of $\Q G$ becomes equal to $\Q G e(G,H,K)$. Hence, if $(H,K)$ is a Shoda pair of $G$, then without using any terminology of the characters, one can simply say that $\Q G e(G, H, K)$ is a simple component of $\Q G$. One has to be little careful that $e(G, H, K)$ is not a primitive central idempotent  but it differs from the same by a rational multiple when $(H,K)$ is a Shoda pair  of $G$.  The case when this rational multiple is $1$, i.e., when $e(G, H, K)$ becomes a primitive central idempotent of $\Q G$ arises with the  following additional constraints on $(H, K)$:
\\ {\rm(i)}  $H \unlhd \C _{G}(\varepsilon(H,K));$ \\ {\rm(ii)}
$\varepsilon(H, K)\varepsilon(H, K)^{g}=0$ $ \forall~ g \in G \setminus \C _{G}(\varepsilon(H,K)),$  where $\varepsilon(H, K)^{g}= g^{-1}\varepsilon(H, K) g$.\vspace{.2cm}\\  The above constraints are sufficient to ensure that $e(G, H, K)$ is a primitive central idempotent of $\Q G$ but are not necessary. A Shoda pair $(H, K)$ of $G$ which satisfy (i) and (ii) above is called \textit{strong Shoda pair} and they were  introduced by  Olivieri, del R{\'{\i}}o and Sim{\'o}n in \cite{OdRS04}. One of the important features of a strong Shoda pair  $(H,K)$ of $G$ (which is proved in Proposition 3.4 of \cite{OdRS04}) is that one can give a precise description of the simple component $\Q G e(G, H, K)$ of $\Q G$ as a matrix algebra over a cyclotomic algebra, which is theoretically known to exist from Brauer Witt theorem. \par  We gave a generalization of strong Shoda pairs in \cite{BK2} and called it generalized strong Shoda pairs.  We say that a Shoda  pair $(H,K)$  of $G$ is  a \textit{generalized strong Shoda pair} of $G$ if there is a chain $H=H_{0}\leq H_{1}\leq \cdots \leq H_{n}=G$ (called strong inductive chain from $H$ to $G$) of subgroups of $G$ such that the following conditions hold for all $ 0 \leq i \leq n-1$: \begin{description} \item [(i)] $H_i \unlhd \operatorname{Cen}_{H_{i+1}}(e_{\mathbb{Q}}(\lambda^{H_i}))$; \item [(ii)] the distinct $H_{i+1}$-conjugates of $e_{\mathbb{Q}}(\lambda^{H_i})$ are mutually orthogonal.\end{description} Here $\lambda$  is a linear character of $H$ with kernel $K$.  Observe that for any $i$, $e_{\mathbb{Q}}(\lambda^{H_i})$ is a rational multiple of $e(H_i, H, K)$ and so their centralizers in $H_{i+1}$ coincide. Thus, one can replace conditions (i) and (ii) above with the following equivalent conditions: \begin{description} \item [(i)] $H_i \unlhd \operatorname{Cen}_{H_{i+1}}(e(H_{i}, H, K))$; \item [(ii)] the distinct $H_{i+1}$-conjugates of $e(H_{i}, H, K)$ are mutually orthogonal.\end{description}  In Theorem 3  of \cite{BK2}, we gave an explicit  description of the structure of the simple component $\Q G e(G, H, K)$ when $(H, K)$ is a generalized strong Shoda pair of $G$.  \par  A group $G$ is called monomial (resp. strongly monomial/ generalized strongly monomial) if every complex irreducible character of $G$ comes from a Shoda pair (resp. strong  Shoda pair/generalized strong Shoda pair) of $G$.  Clearly strongly monomial groups $\subset$ generalized strongly monomial groups $\subset $ monomial groups. The class of monomial groups is well known in the literature  from  the pioneering  work of  Dade, Isaacs, Gunter, Dornhoff and others. While it is  proved in \cite{OdRS04} that all abelian-by-supersolvable groups are strongly monomial,  in \cite{BK4} we provided an extensive list of the groups which are generalized strongly monomial.  In particular, the class of generalized strongly monomial groups include strongly monomial groups, subnormally monomial groups, nilpotent-by-supersolvable monomial groups of odd order and monomial groups with Sylow tower. Moreover, in \cite{BK4}, we revisited Dade's embedding theorem which states that every finite solvable group can be embedded in some monomial group and proved that the embedding is indeed done in some generalized strongly monomial group.
\par Two generalized strong Shoda pairs of $G$ are said to be equivalent if they realize the same primitive central idempotent of $\mathbb{Q}G$. A set of representatives of distinct equivalence classes of  generalized strong Shoda pairs of $G$ is called a \textit{complete and irredundant set of generalized strong Shoda pairs} of $G$.  \par
For generalized strongly monomial groups, some results were obtained in \cite{BK2} related to the study of the unit group of integral group ring. We need to recall them in order to come to the work in this paper. Firstly, let's recall Bass units and generalized Bass units. Given $g \in G$ and $k,m$ positive integers such that $k^{m}\equiv 1\operatorname{mod}|g|$, where $|g|$ is the order of $g$, the following is a unit of $\mathbb{Z}G$: $$u_{k,m}(g)=(1+g+\cdots +g^{k-1})^{m}+\frac{1-k^{m}}{|g|}(1+g+\cdots +g^{|g|-1}).$$ The units of this form are called \textit{Bass units} based on $g$ with parameters $k$ and $m$ and were introduced by Bass \cite{Bass}. When $M$ is a normal subgroup of $G$, then $$u_{k,m}(1-\widehat{M}+g\widehat{M})=1-\widehat{M}+u_{k,m}(g)\widehat{M}$$ is an invertible element of $\mathbb{Z}G(1-\widehat{M})+\mathbb{Z}G\widehat{M}$. As this is an order in $\mathbb{Q}G$, for each element $b=u_{k,m}(1-\widehat{M}+g\widehat{M})$ there is a positive integer $n$ such that $b^{n} \in \mathcal{U}(\mathbb{Z}G)$. Let $n_{b}$ denote the minimal positive integer satisfying this condition. The element $$u_{k,m}(1-\widehat{M}+g\widehat{M})^{n_{b}}=1-\widehat{M}+u_{k,mn_{b}}(g)\widehat{M}$$ is called the \textit{generalized Bass unit} \cite{JOdRV} based on $g$ and $M$ with parameters $k$ and $m$. When $(H, K)$ is a generalized strong Shoda pair of $G$, then we provided, in \cite{BK2}, an iterative process to construct central units of $\mathbb{Z}G$ from the central units of $\mathbb{Z}H$ lying in $\mathcal{Z}(\mathcal{U}(\mathbb{Z}(1-\varepsilon(H,K))+\mathbb{Z}H\varepsilon(H,K)))$. Let $\lambda$ be a linear character of $H$ with kernel $K$ and let $H=H_{0} \leq H_{1} \leq \cdots \leq H_{n}=G$ be a strong inductive chain from $H$ to $G$. For $u \in \mathcal{Z}(\mathcal{U}(\mathbb{Z}H))\cap \mathcal{Z}(\mathcal{U}(\mathbb{Z}(1-\varepsilon(H,K))+\mathbb{Z}H\varepsilon(H,K)))$, put $$z_{0}^{\mathcal{N}}(u)=u$$ and for $0\leq i\leq n-1$, put $$z_{i+1}^{\mathcal{N}}(u)=\prod_{t \in T_{i}}\big(\prod_{c \in C_{i}}z_{i}^{\mathcal{N}}(u)^{c}\big)^{t},$$ where $C_{i}=\operatorname{Cen}_{H_{i+1}}(e_{\mathbb{Q}}(\lambda^{H_{i}}))$ and $T_{i}$ is a right transversal of $C_{i}$ in $H_{i+1}$. Denote the final step of the construction $z_{n}^{\mathcal{N}}(u)$ by $z^{\mathcal{N}}(u)$.  In \cite{BK2}, it is proved that $z^{\mathcal{N}}(u)$ is a central unit of $\mathbb{Z}G$ and certain suitable product of units of this kind, when $u$ is taken in the subgroup generated by generalized Bass units,  generate a large subgroup of $\mathcal{Z}(\mathcal{U}(\mathbb{Z}G))$  in the case when $G$ is a generalized strongly monomial group with a constraint on a complete and irredundant set of Shoda pairs.  Let's see  a terminology (termed as least center property w.r.t. a strong inductive chain) used in the constraint. When $(H,K)$ is a generalized strong Shoda pair of $G$ and $\mathcal{N}:~H=H_{0}\leq H_{1}\leq \cdots \leq H_{n}=G$ is a strong inductive chain from $H$ to $G$, then we proved in \cite{BK2} that, for all $0 \leq i \leq n-1$,  $\Q H_{i+1}e_{\mathbb{Q}}(\lambda^{H_{i+1}})$ is isomorphic to a matrix algebra over the crossed product $\Q H_{i}e_{\mathbb{Q}}(\lambda^{H_{i}}) *^\sigma _{\tau} C_i/H_i$ with some suitably defined action $\sigma$ and twisting $\tau$. In case, the center of  $\Q H_{i+1}e_{\mathbb{Q}}(\lambda^{H_{i+1}})$ is precisely all the elements of the center of $\Q H_{i}e_{\mathbb{Q}}(\lambda^{H_{i}})$ which are kept fixed by the action of $C_i/H_i$, then we say that $(H,K)$  has the least center property w.r.t.  the strong inductive chain $\mathcal{N}$ from $H$ to $G$. We are now ready to  precisely  state the result proved in \cite{BK2}: \begin{theorem*}\rm{(\cite{BK2}, Theorem 4)} Let $G$ be a generalized strongly monomial group and let $\mathcal{S}=\{(H_{i},K_{i})~|~1\leq i \leq n\}$ be a complete and irredundant set of generalized strong Shoda pairs of $G$. For each $i$, let $A_{(H_{i},K_{i})}$ be the subgroup of $\mathcal{Z}(\mathcal{U}(\mathbb{Z}H_{i}))$ generated by the generalized Bass units $b_{i}^{n_{b_{i}}}$, where $b_{i}=u_{k,m}(1-\widehat{H_{i}'}+h\widehat{H_{i}'})$, $h \in H_{i}$, $k$ and $m$ positive integers s.t. $k^{m}\equiv 1\operatorname{mod}|h|$. If $(H_{i},K_{i})$ has the least center property w.r.t. a strong inductive chain $\mathcal{N}_{i}$ from $H_{i}$ to $G$ for all $1\leq i \leq n$, then $$\small\{z^{\mathcal{N}_{1}}(u_{1})\cdots z^{\mathcal{N}_{n}}(u_{n})~|~u_{i}\in A_{(H_{i},K_{i})}\cap \mathcal{Z}(\mathcal{U}(\mathbb{Z}(1-\varepsilon(H_{i},K_{i}))+\mathbb{Z}H_{i}\varepsilon(H_{i},K_{i}))),~1\leq i\leq n\}$$ forms a subgroup of $\mathcal{Z}(\mathcal{U}(\mathbb{Z}G))$ which is contained in the group generated by  generalized Bass units of $\mathbb{Z}G$ and its index in $\mathcal{Z}(\mathcal{U}(\mathbb{Z}G))$ is finite.\end{theorem*}  If $g\in G$ is such that $\langle g \rangle$ is subnormal in $G$, then a construction of central unit of $\mathbb{Z}G$ from the Bass unit based on $g$ is given in \cite{JOdRV}. The same idea was  used in \cite{BK2} to provide a construction of central units of $\mathbb{Z}G$ from central units of $\mathbb{Z}H$, when $H$ is a subnormal subgroup of $G$. Let $\mathcal{N}:H=H_{0}\unlhd H_{1}\unlhd \cdots \unlhd H_{n}=G$ be a subnormal series. For $u \in \mathcal{Z}(\mathcal{U}(\mathbb{Z}H))$, put $c_{0}^{\mathcal{N}}(u)=u$ and for $0\leq i\leq n-1$, set $$c_{i+1}^{\mathcal{N}}(u)=\prod_{t \in T_{i}}c_{i}^{\mathcal{N}}(u)^{t},$$ where $T_{i}$ is a right transversal of $H_{i}$ in $H_{i+1}$. Denote $c_{n}^{\mathcal{N}}(u)$ by $c^{\mathcal{N}}(u)$. It was proved that $c^{\mathcal{N}}(u)$ is well defined and is a central unit of $\mathbb{Z}G$.  Such type of central units  have been proved to generate a large subgroup of  $\mathcal{Z}(\mathcal{U}(\mathbb{Z}G))$ when $G \in \mathcal{C}$ but with a constraint on a complete and irredundant set of Shoda pairs as stated in the following:  \begin{theorem*}\rm{(\cite{BK2}, Theorem 5)} Let $G \in \mathcal{C}$ be such that every cyclic subgroup of order not a divisor of $4$ or $6$ is subnormal in $G$.  Let $S$ be a complete and irredundant set of generalized strong Shoda pairs of $G$. If each $(H,K) \in S$ has the least center property w.r.t. some strong inductive chain from $H$ to $G$, then,
 $$ \langle c^{\mathcal{N}_{g}}(b_{g})~|~b_{g}~{\rm is~a~Bass~unit~based~on}~g \in G ~{\rm with~order}~{\rm not~divisble ~by}~4~{\rm or }~6\rangle $$ is a subgroup of the group generated by Bass units of $\mathbb{Z}G$ which is of finite index in $\mathcal{Z}(\mathcal{U}(\mathbb{Z}G)),$ where $\mathcal{N}_{g}$ is a fixed subnormal series from the cyclic subgroup $\langle g \rangle$ to $G$. \end{theorem*} In the next section, we will show that every generalized strong Shoda pair of $G$ has the so called least center property w.r.t. any strong inductive chain and hence the extra condition imposed on generalized strong Shoda pairs in both of the above theorems is no longer needed.
\section{Large subgroups of $\mathcal{Z}(\mathcal{U}(\mathbb{Z}G))$} The following are the two main results proved in this section. \begin{theorem}\label{t1}Let $G$ be a subgroup closed monomial group such that every cyclic subgroup of order not a divisor of $4$ or $6$ is subnormal in $G$. Then  the group generated by Bass units of $\mathbb{Z}G$ contains a large subgroup of  $\mathcal{Z}(\mathcal{U}(\mathbb{Z}G))$.\end{theorem} \begin{theorem}\label{t2}Let $G$ be a  generalized  strongly monomial group. Then the group generated by generalized Bass units of $\mathbb{Z}G$ contains a large subgroup of  $\mathcal{Z}(\mathcal{U}(\mathbb{Z}G))$. \end{theorem}\noindent   A careful examination of the proof of ($\cite{BK2}$, Theorem 5) reveal that its  arguments work  well when $G$ in $\mathcal{C}$ is replaced by any generalized strongly monomial group all of whose subgroups are generalized strongly monomial.  In a recent work done in \cite{BK4}, we proved that subgroup closed monomial groups are generalized strongly monomial. Hence, one obtains  Theorem \ref{t1} once we show that each generalized strong Shoda pair of $G$ has the so called least center property w.r.t. any strong inductive chain.  In view of (\cite{BK2}, Theorem 4), the same is the requirement to prove Theorem \ref{t2}. Hence both the above theorems will be proved once we prove the following: \begin{prop}\label{p1} Let $G$ be a finite group. Let $(H,K) $  be a generalized strong Shoda pair of $G$ and let $H=H_{0}\leq H_{1}\leq \cdots \leq H_{n}=G$ be a strong inductive chain from $H$ to $G$.  Let $\lambda$  be a linear character of $H$ with kernel $K$.  Then  the following hold for all $0 \leq i \leq n-1$: \begin{enumerate} \item   $\mathbb{Q}C_{i}e_{\mathbb{Q}}(\lambda^{H_{i}}) \cong  \mathbb{Q}H_{i}e_{\mathbb{Q}}(\lambda^{H_{i}})\ast^{\sigma_{H_{i}}}_{\tau_{H_{i}}}C_{i}/H_{i},$ where $\sigma_{H_{i}}:C_{i}/H_{i}\rightarrow \operatorname{Aut}(\mathbb{Q}H_{i}e_{\mathbb{Q}}(\lambda^{H_{i}}))$ maps $x$ to the conjugation automorphism $(\sigma_{H_{i}})_{x}$ on $\mathbb{Q}H_{i}e_{\mathbb{Q}}(\lambda^{H_{i}})$ induced by $\overline{x}$ and  $\tau_{H_{i}}:
		C_{i}/H_{i}\times C_{i}/H_{i}\rightarrow \mathcal{U}(\mathbb{Q}H_{i}e_{\mathbb{Q}}(\lambda^{H_{i}}))$ is given by $\tau_{H_{i}}(x,y) = \overline{x}.\overline{y}.\overline{xy}^{-1}e_{\mathbb{Q}}(\lambda^{H_{i}}).$ Here  $C_i = \operatorname{Cen}_{H_{i+1}}(e_{\mathbb{Q}}(\lambda^{H_{i}}))$ and for each $x \in C_{i}/H_{i}$, $\overline{x} \in C_i$ is a fixed inverse image of  $x \in C_i/H_i$ under the natural map $C_{i} \rightarrow C_{i}/H_{i}$. Furthermore, $\sigma_{H_{i}}$   is not an  inner automorphism;
		\item $C_i/H_i$ acts faithfully on the center of  $\mathbb{Q}H_{i}e_{\mathbb{Q}}(\lambda^{H_{i}})$ by conjugation.
	\end{enumerate}\end{prop}\noindent {\bf Proof.} (i) Observe that  $\mathbb{Q}C_{i}e_{\mathbb{Q}}(\lambda^{H_{i}}) $ is isomorphic to the crossed product $\mathbb{Q}H_{i}e_{\mathbb{Q}}(\lambda^{H_{i}})$\linebreak $\ast^{\sigma_{H_{i}}}_{\tau_{H_{i}}}C_{i}/H_{i}$ by taking $R= \mathbb{Q}H_{i}e_{\mathbb{Q}}(\lambda^{H_{i}}), N=H_{i}$ and $G=C_{i}$ in Lemma 2.6.2 of \cite{JdR}. Now we show that  $(\sigma_{H_{i}})_{x}$ is not an inner automorphism of $\mathbb{Q}H_{i}e_{\mathbb{Q}}(\lambda^{H_{i}})$  for all $x \in C_i/H_i$ and $ 0\leq i \leq n-1$. Suppose that for some $0 \leq i \leq n-1$, there exists $x \in C_{i}/H_{i}$ such that $(\sigma_{H_{i}})_{x}$ is an inner automorphism of $\mathbb{Q}H_{i}e_{\mathbb{Q}}(\lambda^{H_{i}})$, i.e., there is a unit $u \in \mathbb{Q}H_{i}e_{\mathbb{Q}}(\lambda^{H_{i}})$ such that \begin{equation}\label{e1} \overline{x}\alpha\overline{x}^{-1} = u\alpha u^{-1} ~~~~~\forall~\alpha \in \mathbb{Q}H_{i}e_{\mathbb{Q}}(\lambda^{H_{i}}).\end{equation} As the above eqn holds for all
$\alpha \in \mathbb{Q}H_{i}e_{\mathbb{Q}}(\lambda^{H_{i}})$, in particular, it holds for $\alpha = he_{\mathbb{Q}}(\lambda^{H_{i}})$, where  $h \in H_{i}$. Hence \begin{equation}\label{e2} \overline{x}he_{\mathbb{Q}}(\lambda^{H_{i}})\overline{x}^{-1} = u h e_{\mathbb{Q}}(\lambda^{H_{i}}) u^{-1}  ~~~~~\forall ~h \in H_{i}. \end{equation} Let $\rho_{H_{i}}$ be a representation of $H_{i}$ affording the character $\lambda^{H_{i}}$. Extending $\rho_{H_{i}}$ linearly on $\mathbb{Q}H_{i}e_{\mathbb{Q}}(\lambda^{H_{i}})$ and applying on eqn \ref{e2}, we get $$\rho_{H_{i}}(\overline{x}h\overline{x}^{-1})= \rho_{H_{i}}(u)\rho_{H_{i}}(h)\rho_{H_{i}}(u)^{-1}~~~\forall~h \in H_{i}.$$ By taking trace, we have $$\lambda^{H_{i}}(\overline{x}h\overline{x}^{-1})= \lambda^{H_{i}}(h) ~~~~~\forall~h\in H_{i}.$$  Thus we have obtained  that $\overline{x}$ stabilizes $\lambda^{H_{i}}$ in $C_{i}$. Since $H_{i}$ is normal in $C_{i}$ and $\lambda^{H_{i}}$ is irreducible, it follows from Mackey's irreducibility criterion that $\overline{x} \in H_{i}$, i.e., $x$ is identity and hence (i) is proved. \\ (ii) Suppose $x \in C_{i}/H_{i}$ such that $\overline{x}\alpha \overline{x}^{-1} =\alpha $ for all $\alpha \in \mathcal{Z}(\mathbb{Q}H_{i}e_{\mathbb{Q}}(\lambda^{H_{i}}))$. Observe that $\sum_{g \in T}h^{g}e_{\mathbb{Q}}(\lambda^{H_{i}}) \in \mathcal{Z}(\mathbb{Q}H_{i}e_{\mathbb{Q}}(\lambda^{H_{i}}))$, where $h \in H_{i}$ and  $T$ is a right transversal of $\operatorname{Cen}_{H_{i}}(h)$ in $H_{i}$. Let $\rho_{i}$ be a representation afforded by $\lambda^{H_{i}}$. Extending $\rho_{i}$ linearly on $\mathbb{Q}H_{i}$, we have $\rho_{i}(\overline{x}(\sum_{g \in T}h^{g})\overline{x}^{-1})= \rho_{i}(\sum_{g \in T}h^{g})$. Hence $(\lambda^{H_{i}})^{\overline{x}}(h)=\lambda^{H_{i}}(h)$ which yields that $\overline{x}$ stabilizes $\lambda^{H_{i}}$ in $C_{i}$. Consequently $\overline{x} \in H_{i}$ and thus $x$ is identity.  This completes the proof. \qed \begin{cor} Let $G$ be a finite group and let $(H,K)$ be a generalized strong Shoda pair  of $G$. Then  $(H,K)$ has the least center property w.r.t. any strong inductive chain from $H$ to $G$.  \end{cor} \noindent {\bf Proof.} Let  $H=H_{0}\leq H_{1}\leq \cdots \leq H_{n}=G$  be a strong inductive chain from $H$ to $G$ and $\lambda$ be a linear character of $H$ with kernel $K$. Let $C_i = \operatorname{Cen}_{H_{i+1}}(e_{\mathbb{Q}}(\lambda^{H_{i}}))$. In Proposition 2 of \cite{BK2}, we have proved that for all $0 \leq i \leq n-1$,  $\mathbb{Q}H_{i+1}e_{\mathbb{Q}}(\lambda^{H_{i+1}})$ is isomorphic to a matrix algebra over  $\mathbb{Q}C_{i}e_{\mathbb{Q}}(\lambda^{H_{i}}) $. Hence $\mathcal{Z}(\mathbb{Q}H_{i+1}e_{\mathbb{Q}}(\lambda^{H_{i+1}})) \cong  \mathcal{Z}(\mathbb{Q}C_{i}e_{\mathbb{Q}}(\lambda^{H_{i}}))$. Proposition \ref{p1} gives that $\mathbb{Q}C_{i}e_{\mathbb{Q}}(\lambda^{H_{i}}) \cong  \mathbb{Q}H_{i}e_{\mathbb{Q}}(\lambda^{H_{i}})\ast^{\sigma_{H_{i}}}_{\tau_{H_{i}}}C_{i}/H_{i}$, where $ \sigma_{H_{i}}$ is not an inner automorphism.   Hence, using Lemma  2.6.1 of \cite{JdR}, we have that $\mathcal{Z}(\mathbb{Q}C_{i}e_{\mathbb{Q}}(\lambda^{H_{i}}) )\cong  \mathcal{Z}(\mathbb{Q}H_{i}e_{\mathbb{Q}}(\lambda^{H_{i}})) ^{C_i/H_i}$. This yields the desired result.\qed

 \section{The rank of $\mathcal{Z}(\mathcal{U}(\mathbb{Z}G))$}
\begin{theorem}\label{t3} Let $G$ be a generalized strongly monomial group and let $\mathcal{S}$ be a complete and irredundant set of generalized strong Shoda pairs of $G$.  For each $(H,K) \in \mathcal{S}$, fix a strong inductive chain: $H=H_{0}\leq H_{1}\cdots \leq H_{n_{(H,K)}}=G$ from $H$ to $G$. Then  the rank of $\mathcal{Z}(\mathcal{U}(\mathbb{Z}G))$ equals $$ \sum_{(H,K)\in \mathcal{S}}\frac{\phi([H:K])}{k_{(H,K)}[C_{0}:H_{0}]  \cdots [C_{n_{(H,K)}-1}:H_{n_{(H,K)}-1}]}-1,$$ where $C_{i}$ is the centralizer of $e(H_{i},H,K)$ in $H_{i+1}$ and $$ k_{(H,K)}:=\left\{\begin{array}{ll}1, & \hbox{$\mathcal{Z}(\mathbb{Q}Ge(G,H,K))$ is totally real;} \\2, & \hbox{otherwise}.\end{array}\right. $$  \end{theorem}\noindent{\bf Proof.}   Let $(H,K) \in \mathcal{S}$ and let $\lambda$ be a linear character of $H$ with kernel $K$. In view of Proposition \ref{p1}, we have that, for all $0 \leq i \leq  n_{(H,K)}-1$,  \begin{equation}\label{e3}\mathcal{Z}(\mathbb{Q}H_{i+1} e_{\mathbb{Q}}(\lambda^{H_{i+1}}))\cong  \mathcal{Z}(\mathbb{Q}H_{i} e_{\mathbb{Q}}(\lambda^{H_{i}})) ^{C_{i}/H_{i}}.\end{equation}  Recursively using eqn \ref{e3} yields that $$  \mathcal{Z}(\mathbb{Q}G e_{\mathbb{Q}}(\lambda^{G}))=\mathcal{Z}(\mathbb{Q}Ge(G,H,K))\cong   ((\mathbb{Q}(\xi_{[H:K]})^{C_{0}/H_{0}})^{C_{1}/H_{1}}\cdots)^{C_{n-1}/H_{n-1}}.$$ The  Wedderburn decomposition $\mathbb{Q}G \cong \oplus_{(H,K) \in \mathcal{S}} \mathbb{Q}G e(G,H,K)$  thus implies that  $$\mathcal{Z}(\mathbb{Q}G)\cong \displaystyle{\oplus_{(H,K) \in \mathcal{S}}} ((\mathbb{Q}(\xi_{[H:K]})^{C_{0}/H_{0}})^{C_{1}/H_{1}}\cdots)^{C_{n_{(H,K)}-1}/H_{n_{(H,K)}-1}}.$$   For notational convenience, let's denote $((\mathbb{Q}(\xi_{[H:K]})^{C_{0}/H_{0}})^{C_{1}/H_{1}}\cdots)^{C_{n_{(H,K)}-1}/H_{n_{(H,K)}-1}}$ by $\mathbb{Q}_{(H,K)}$ and $((\mathbb{Z}[\xi_{[H:K]}]^{C_{0}/H_{0}})^{C_{1}/H_{1}}\cdots)^{C_{n_{(H,K)}-1}/H_{n_{(H,K)}-1}}$ by $\mathbb{Z}_{(H,K)}$. Observe that $\mathbb{Z}_{(H,K)}$ is the ring of integers of  $\mathbb{Q}_{(H,K)}$ and  it is  thus the unique maximal order of $\mathbb{Q}_{(H,K)}.$ Since any maximal order in a semisimple algebra is a direct sum of maximal orders in its simple constituents,  it follows that $\oplus_{(H,K) \in \mathcal{S}}\mathbb{Z}_{(H,K)}$ is the unique maximal order of $\oplus_{(H,K) \in \mathcal{S}}\mathbb{Q}_{(H,K)}\cong \mathcal{Z}(\mathbb{Q}G)$.   Now $\mathcal{Z}(\mathbb{Z}G)$ being a $\mathbb{Z}$-order in $\mathcal{Z}(\mathbb{Q}G)$, we obtain that $\mathcal{Z}(\mathbb{Z}G)$   is contained in $\oplus_{(H,K) \in \mathcal{S}}\mathbb{Z}_{(H,K)}.$ Consequently Lemma 4.6.9 of \cite{JdR} yields that  $[\mathcal{U}(\oplus_{(H,K)}\mathbb{Z}_{(H,K)}): \mathcal{U}(\mathcal{Z}(\mathbb{Z}G))] < \infty$. Also, using  Dirichlet's unit theorem, $\mathcal{U}(\mathbb{Z}_{(H,K)})$ is a finitely generated abelian group. Hence the rank of   $\mathcal{U}(\mathcal{Z}(\mathbb{Z}G))$ is equal to the the sum of ranks of  $\mathcal{U}(\mathbb{Z}_{(H,K)})$, where the sum runs over  $(H,K) \in \mathcal{S}.$ As $\mathbb{Q}_{(H,K)}$  is a Galois extension over $\mathbb{Q}$, we have that $\mathbb{Q}_{(H,K)}$ is  either totally real or totally imaginary. Hence, by Dirichlet's unit theorem, the rank of  $\mathcal{U}(\mathbb{Z}_{(H,K)})$  is equal to  $\frac{[\mathbb{Q}_{(H,K)}:\mathbb{Q}]}{k_{(H,K)}}-1$, where $k_{(H,K)}=1$ if $\mathbb{Q}_{(H,K)}$ is totally real and $k_{(H,K)}=2$ if $\mathbb{Q}_{(H,K)}$ is totally imaginary.  The faithful action of $C_{i+1}/H_{i+1}$ on  $((\mathbb{Q}(\xi_{[H:K]})^{C_{0}/H_{0}})^{C_{1}/H_{1}}\cdots)^{C_{i}/H_{i}}$ by conjugation for all $0 \leq i \leq n-1$ yields that $[\mathbb{Q}_{(H,K)}:\mathbb{Q}] = \frac{\phi([H:K])}{[C_{0}:H_{0}]  \cdots [C_{n_{(H,K)}-1}:H_{n_{(H,K)}-1}]}$. This finishes the proof of the theorem. \qed \vspace{.2cm}\\\noindent  Observe that for any Shoda pair $(H,K)$ of a group $G$, $\mathcal{Z}(\mathbb{Q}Ge(G,H,K))$  is totally real if and only if $\lambda^G$ is a real character, where $\lambda$ is a linear character of $H$ with kernel $K$. Hence, the following corollary arises: \begin{cor}\label{t4} If $G$ is a real group which is generalized strongly monomial, then the rank of $\mathcal{Z}(\mathcal{U}(\mathbb{Z}G))$ equals $$\sum_{(H,K)\in\mathcal{S}}\frac{\phi([H:K])}{[C_{0}:H_{0}] \cdots [C_{n_{(H,K)}-1}:H_{n_{(H,K)}-1}]}-1,$$   where the notations are as in the the above theorem.\end{cor}\vspace{.2cm} \noindent A necessary and sufficient condition for a group $G$ to have a non trivial real valued irreducible character is that the order of the group $G$ is even. This gives rise to the following: \begin{cor}\label{t5} If $G$ is a generalized strongly monomial group of odd order,  then the rank of $\mathcal{Z}(\mathcal{U}(\mathbb{Z}G))$ equals $$ \sum_{(H,K)\in S} \frac{\phi([H:K])}{2[C_{0}:H_{0}] \cdots [C_{n_{(H,K)}-1}:H_{n_{(H,K)}-1}]}-1,$$  where  the notations are as explained above. \end{cor}    \noindent Now, we illustrate the above theory with an example: \vspace{.25cm}\\ \textbf{Example}:  We will compute the rank of $\mathcal{Z}(\mathcal{U}(\mathbb{Z}G))$, where $G$=SmallGroup(1000,86) in  GAP library, the smallest group which is generalized strongly monomial but not strongly monomial (see \cite{BK1}). This group is a semidirect product of the extraspecial $5$-group of order $5^{3}$ with the cyclic group of order 8. It is generated by $x_{i},~1\leq i\leq 6$, with the following defining relations:
\begin{quote}  $x_{1}^{2}x_{2}^{-1}$=$x_{2}^{2}x_{3}^{-1}$=$x_{4}^{5}$=$x_{3}^{2}$=$x_{5}^{5}$=$x_{6}^{5}$=$1,$\\
               $[x_{2},x_{1}]$=$[x_{3},x_{1}]$=$[x_{3},x_{2}]$=$[x_{6},x_{3}]$=$[x_{6},x_{4}]$=$[x_{6},x_{5}]$=$1,$\\
                $[x_{5},x_{4}]$=$x_{6}$, $[x_{5},x_{1}]$=$x_{4},$\\
                $[x_{6},x_{1}]$=$x_{6}^{2}$, $[x_{4},x_{2}]$=$x_{4}x_{6}^{2}$, $[x_{6},x_{2}]$=$x_{6}^{3},$\\
                 $[x_{5},x_{2}]$=$x_{5}x_{6}^{2}$, $[x_{5},x_{3}]$=$x_{5}^{3}x_{6}^{2}$, $[x_{4},x_{3}]$=$x_{4}^{3}x_{6}^{2}$, $[x_{4},x_{1}]$=$x_{4}^{-2}x_{5}^{3}x_{6}^{4}.$\end{quote}
                 In section 7 of \cite{BK1}, we proved that the following set $S$ of Shoda pairs of $G$ form a complete and irredundant set of Shoda pairs: $\{(G,G)$, $(G,\langle x_{2},x_{3},x_{4},x_{5},x_{6}\rangle)$, $(G,\langle x_{3},x_{4},x_{5},x_{6}\rangle),$ $(G,\langle x_{4},x_{5},x_{6} \rangle),$ $(\langle x_{4},x_{5},x_{6} \rangle,\langle x_{4},x_{6}\rangle),$ $(\langle x_{4},x_{5},x_{6} \rangle,\langle x_{5},x_{6}\rangle)$,\linebreak $(\langle x_{4},x_{5},x_{6} \rangle,\langle x_{4}^{-1}x_{5},x_{6}\rangle),$ $(\langle x_{5},x_{6},x_{3}x_{4}^{2}x_{6}\rangle, \langle x_{3}x_{4}^{2}x_{6}^{3},x_{5}\rangle),$ $(\langle x_{5},x_{6},x_{3}x_{4}^{2}\rangle,\langle x_{5}\rangle)\}$.  Observe that, in fact $\mathcal{S}$ is a complete and irredundant set of generalized strong Shoda pairs of $G$ and all pairs except the last two are strong Shoda pairs of $G$. In view of Theorem \ref{t3}, to determine the rank of $\mathcal{Z}(\mathcal{U}(\mathbb{Z}G))$, for each generalized strong Shoda pair $(H,K)\in S$, we need to determine a strong inductive chain $\mathcal{N}:H=H_{0} \leq H_{1}\leq \cdots \leq H_{n}=G$  from $H$ to $G$, $[C_{i}:H_{i}]$ for all $0 \leq i \leq n-1$, $[H:K]$ and $k_{(H,K)}$. Notice that if $(H,K)$ is a strong Shoda pair of $G$ then $"H\leq G"$ is a strong inductive chain. \vspace{.2cm}\\ $\underline{(H,K)=(G,G)}$ \vspace{.15cm}\\ $[H:K]=1$, $[C_{0}:H_{0}]=1$ and $k_{(G,G)}=1$. \vspace{.2cm}\\ $\underline{(H,K)=(G,\langle x_{2},x_{3},x_{4},x_{5},x_{6}\rangle)}$ \vspace{.15cm}\\   Here $[H:K]=2$ and $[C_{0}:H_{0}]=1$. As the index of $\langle x_{2},x_{3},x_{4},x_{5},x_{6}\rangle$ in $G$ is 2, the irreducible character of $G$ with kernel
                 $\langle x_{2},x_{3},x_{4},x_{5},x_{6}\rangle$ is real. Hence the simple component associated to this character is totally real and thus $k_{(G,\langle x_{2},x_{3},x_{4},x_{5},x_{6}\rangle)}=1$. \vspace{.2cm}\\ $\underline{(H,K)=(G,\langle x_{3},x_{4},x_{5},x_{6}\rangle)}$ \vspace{.15cm}\\ $[H:K]=4$ and $[C_{0}:H_{0}]=1$. The index of $\langle x_{3},x_{4},x_{5},x_{6}\rangle$ in $G$ being 4 yields that $k_{(G,\langle x_{3},x_{4},x_{5},x_{6}\rangle)}=2.$ \vspace{.2cm}\\ $\underline{(H,K)=(G,\langle x_{4},x_{5},x_{6}\rangle)}$ \vspace{.15cm}\\ $[H:K]=8$, $[C_{0}:H_{0}]=1$ and $k_{(G,\langle x_{4},x_{5},x_{6}\rangle)}=2$, as the index of $\langle x_{4},x_{5},x_{6}\rangle$ in $G$ is 8. \vspace{.2cm}\\ $\underline{(H,K)=(\langle x_{4},x_{5},x_{6}\rangle, \langle x_{4},x_{6}\rangle)}$ \vspace{.15cm}\\ In this case $[H:K]=5$ and $[C_{0}:H_{0}]=4$. It can be checked that monomial character induced from a linear character of subgroup $\langle x_{4},x_{5},x_{6} \rangle$ with kernel $\langle x_{4},x_{6}\rangle$ is real. Therefore, the simple component associated to this monomial real character is totally real. Hence $k_{(\langle x_{4},x_{5},x_{6}\rangle,\langle x_{4},x_{6}\rangle)}=1$. \vspace{.2cm}\\ $\underline{(H,K)=(\langle x_{4},x_{5},x_{6}\rangle, \langle x_{5},x_{6}\rangle)}$ \vspace{.15cm}\\ $ [H:K]=5,$ $[C_{0}:H_{0}]=4$ and $k_{(\langle x_{4},x_{5},x_{6}\rangle,\langle x_{5},x_{6}\rangle)}=1$. \vspace{.2cm}\\ $\underline{(H,K)=(\langle x_{4},x_{5},x_{6}\rangle, \langle x_{4}^{-1}x_{5},x_{6}\rangle)}$ \vspace{.15cm}\\ This is also a strong Shoda pair of $G$ with a strong inductive chain from $H$ to $G$ as $H\leq G, [H:K]=5,$ $[C_{0}:H_{0}]=4$ and $k_{(\langle x_{4},x_{5},x_{6}\rangle,\langle x_{4}^{-1}x_{5},x_{6}\rangle)}=1$.
                 \vspace{.2cm}\\ $\underline{(H,K)=(\langle x_{5},x_{6},x_{3}x_{4}^{2}\rangle, \langle x_{3}x_{4}^{2}x_{6}^{3},x_{5}\rangle)}$   \vspace{.15cm}\\ This is a generalized strong Shoda pair of $G$ but not a strong Shoda pair. In this case, $H=
                 \langle x_{5},x_{6},x_{3}x_{4}^{2}\rangle \leq \langle x_{5},x_{6},x_{3}x_{4}^{2}\rangle \leq \langle x_{3},x_{4},x_{5},x_{6} \rangle \leq G$ can be taken as a strong inductive chain. Note that $[H:K]=5$, $[C_{0}:H_{0}]=1,~[C_{1}:H_{1}]=1$ and $[C_{2}:H_{2}]=4.$ With GAP calculations, one can see that $k_{(\langle x_{5},x_{6},x_{3}x_{4}^{2}x_{6}\rangle, \langle x_{3}x_{4}^{2}x_{6}^{3},x_{5}\rangle)}=1$.\vspace{.2cm}\\ $\underline{(H,K)=(\langle x_{5},x_{6},x_{3}x_{4}^{2}\rangle, \langle x_{5}\rangle)}$ \vspace{.15cm}\\  This is also a generalized strong Shoda pair of $G$ but not a strong Shoda pair.  A strong inductive chain from $H$ to $G$ is $H=
                 \langle x_{5},x_{6},x_{3}x_{4}^{2}\rangle \leq \langle x_{5},x_{6},x_{3}x_{4}^{2}\rangle \leq \langle x_{3},x_{4},x_{5},x_{6} \rangle \leq G$, $[H:K]=10$, $[C_{0}:H_{0}]=1,~[C_{1}:H_{1}]=1,$ $[C_{2}:H_{2}]=4$ and $k_{(\langle x_{5},x_{6},x_{3}x_{4}^{2}x_{6}\rangle, \langle x_{5}\rangle)}=1$.    \vspace{.2cm}\\ These computations when substituted in Theorem \ref{t3} yield that the rank of $\mathcal{Z}(\mathcal{U}(\mathbb{Z}G))$ is 1.              \bibliographystyle{amsplain}
\bibliography{BK5}
 \end{document}